\documentclass[12pt,letterpaper]{article}
\usepackage[hmargin=1in,vmargin=1in]{geometry}

\usepackage{amsmath}
\usepackage{amsfonts}
\usepackage{amsthm}
\usepackage{amssymb}
\usepackage{array}
\usepackage{caption}
\usepackage{color}
\usepackage{enumitem}
\usepackage{epic,eepic}
\usepackage{graphics}
\usepackage{graphicx}
\usepackage{hyperref}
\usepackage{mathtools}
\usepackage{rotating}
\usepackage{verbatim}

\begin{document}

\title{Conway's Wizards}
\author{Tanya Khovanova\\MIT}
\maketitle

\begin{abstract}
I present and discuss a puzzle about wizards invented by John H. Conway.
\end{abstract}

\section{Introduction}
\label{sec:intro}

This article is about a puzzle invented by John H. Conway. Conway emailed it to me in June 2009. I wrote about it on my blog \cite{TK1}. I also invented and posted a simplified version \cite{TK2} and a generalized
version \cite{TK3}. But I never published a solution or a discussion. Now the time has come to do so.

The puzzle is introduced in Section~\ref{sec:original}. In the next Section~\ref{sec:discussion} the problem is discussed and the solution is provided. Section~\ref{sec:simplified} presents a simplified variation of Conway's puzzle. Section~\ref{sec:generalized} produces a more advanced generalization of the puzzle.

\section{The Original Puzzle}\label{sec:original}

John Conway sent me a puzzle about wizards, which he invented in the sixties. Here it is:

\begin{quote}
Last night I sat behind two wizards on a bus, and overheard the following:

A: ``I have a positive integral number of children, whose ages are positive integers, the sum of which is the number of this bus, while the product is my own age." \newline
B: ``How interesting! Perhaps if you told me your age and the number of your children, I could work out their individual ages?" \newline
A: ``No." \newline
B: ``Aha! AT LAST I know how old you are!" \newline

Now what was the number of the bus?
\end{quote}

\section{A Discussion and a Solution}\label{sec:discussion}

This is an incredible puzzle. As soon as the number of the bus is known, A's age is revealed. So the task is to figure out A's age with much less information than wizard B has, who clearly has barely enough information.

This is also an under-appreciated puzzle. It is more interesting than it might seem. When someone announces the answer, it is not clear if they solved it completely. No wonder Conway looked at me with suspicion when the day after receiving the puzzle, I told him the number of the bus.

I am a spoiler. I will spoil you and will tell you how to solve this puzzle, but not right now. First, I want to apologize for the seeming rudeness of wizard A. When he says, ``No", he does not mean, ``I do not want to tell you my age and the number of my children." ``No" is the answer to the previous question and it means, ``My age and the number of my children are not enough to deduce their individual ages."

Here is my real question. Why were the wizards riding the 
bus? Could they have been on a trolley? No, they could not have. Conway designs his problems very carefully. In this puzzle he implicitly introduced a good notation. I will follow his suggestion and denote the \textbf{a}ge of the wizard, the number of the \textbf{b}us, and the number of the \textbf{c}hildren by $a$, $b$, and $c$.

\subsection{Examples}

Let us start with some examples. Suppose $b=5$. Here is the list of possible ages for children, and the corresponding age of the wizard and the number of children:
\begin{itemize}
\item 1, 1, 1, 1, and 1. Also, $a=1$ and $c=5$.
\item 1, 1, 1, and 2. Also, $a=2$ and $c=4$.
\item 1, 1, and 3. Also, $a=3$ and $c=3$.
\item 1, 2, and 2. Also, $a=4$ and $c=3$.
\item 1 and 4. Also, $a=4$ and $c=2$.
\item 5. Also, $a=5$ and $c=1$.
\end{itemize}

For the purpose of the theoretical discussion let us ignore the ridiculously small age of the wizard. We see that if the wizard is not four years old, his age is sufficient to determine the number and the ages of the children. In any case, the wizard's age and the number of his children uniquely determine their ages. So wizard A could not have said ``No." From this example we can conclude that the bus number is not 5. Checking other small cases we can assure ourselves that the bus number is greater than 5: $b > 5$.

We can keep increasing the bus number by 1 until we find the bus number for which the wizard A could have said ``No." One might think this gives us the answer. But this puzzle is more interesting than it seems. We will get the right answer this way, and produce the correct bus number and the age of the wizard, but we will miss a big part of the story.

Let us look at large bus numbers. Suppose that $b=21$. The wizard A can say ``No" in this case. Indeed, if he were 96 and had three children, their ages are not uniquely defined. They could be 1, 8 and 12, or alternatively, 2, 3 and 16. But the important question is whether wizard B can now determine that A is 96. Not really. It might be that wizard A is 240 and has three children: 4, 5, and 12, or alternatively, 3, 8, and 10. So 21 cannot be the bus number, but for a completely different reason. Wizard B could not have figured out the age of wizard A.

Does this always happen with large numbers? Let us try 22. If we think about it, we do not even need to do new calculations. The wizard's age is ambiguous again. He can be 96 and have four children of ages 1, 1, 8, and 12, or 1, 2, 3, and 16. Similarly, he can be 240 and have four children 1, 4, 5, and 12, or alternatively, 1, 3, 8, and 10. Do you see what happened? We can add one more child of age one. This way the number of children increases by 1, the sum increases by 1, and the product stays the same. That means if bus $b$ had two possible ages for wizard B, then the same two ages would work for bus $b+1$. Therefore, we do not need to check large numbers over 21: $b < 21$.

Now our problem becomes finite. We just need to try all the possibilities up to 21 to find the answer. It is lucky for Conway that the bus number is uniquely defined in this way. It could have been that as soon as the wizard A said ``No," the age would not be uniquely defined. Or, it could have been that there were several bus numbers for which wizard B could have guessed A's age.

Now to the answer. The bus number is 12. The only age for which the bus number and the number of children do not define the ages of children is 48. The children ages could be 2, 2, 2, and 6, or on the other hand, 1, 3, 4, and 4.

To complete the solution we need to check that all the bus numbers below 12 will not permit wizard A to say no, and bus numbers above 12 will not define the age of wizard A uniquely. I will leave the first task to the reader. I want to show that bus numbers starting from 13 do not work. Indeed, if the bus number is 13, wizard A can say ``no" if he is 48 and has five children, aged either 1, 2, 2, 2, and 6, or 1, 1, 3, 4, and 4. The second possible age for wizard A is 36, with three children aged either 2, 2, 9, or 1, 6, 6. If the bus number is 13 or greater, wizard B cannot figure out the age of wizard A.

\subsection{The number of children}

Here I would like to reflect on the number of children. 

It is obvious from the start that wizard A has more than two children. If he had one child then his/her age would be both the number of the bus and the same as the father's age. While it is unrealistic, in mathematics strange things can happen. The important part is that if wizard A has one child he could not have said ``no." The same is true for two children: their age distribution is uniquely defined by the sum and the product of their ages.

So wizard A has to have at least three children. One common mistake is to assume that wizard A has to have exactly three children. People who make this mistake might well discover the beauty of this puzzle, but they'll miss the bus, so to speak.

They will find that the first case when wizard A says ``no" and has three children, is for bus 13, when wizard A's age is 36 and his children's ages are either 2, 2, 9, or 1, 6, 6. The next bus number, 14, produces two potential ages for wizard A. First, he could be 40 and his children could be aged either 2, 2, 10, or 1, 5, 8. Second, he could be 72 and his children could be aged either 3, 3, 8, or 2, 6, 6. So if someone tells you that the answer is 13, you can safely guess that that person assumed that wizard A has three children.

\section{Simplified Wizards}\label{sec:simplified}

Now that we have discussed the solution to Conway's puzzle, I would like to give you a simpler puzzle that you have to solve on your own. For the sake of continuity, I suppressed reality: assume ages do not need to be realistic.

\begin{quote}
A: ``I have a positive integral number of children, whose ages are positive integers, the sum of which is the number of this bus, while the product is my own age."\newline
B: ``How interesting! Perhaps if you told me the number of your children, I could work out their individual ages?"
A: ``No."
B: ``Aha! AT LAST I know how many children you have!"

What is the number of the bus? While you are at it, what is the age of the wizard?
\end{quote}

Surprisingly, the answer can again be uniquely determined.

\section{Generalized Wizards}\label{sec:generalized}

Now I want to give you a more difficult puzzle:

\begin{quote}
Last night I sat behind two wizards on a bus, and overheard the following:

A: ``I have a positive integral number of children, whose ages are positive integers, the sum of which is the number of this bus, while the product is my own age. Also, the sum of the squares of their ages is the number of dolls in my collection."\newline
B: ``How interesting! Perhaps if you told me your age, the number of your children, and the number of dolls, I could work out your children’s individual ages."\newline
A: ``No."\newline
B: ``Aha! AT LAST I know how old you are!"

Now, what was the number of the bus?
\end{quote}

By the way, by now you should be able to figure out why I prefer ``dolls" over ``baseball cards."

Although I kept the focus of this puzzle on the bus and on the age of the wizard for the sake of continuity with Conway's original puzzle, I had to sacrifice realism. For this puzzle, you need to have an open mind. In Conway's original puzzle you do not need to assume that wizard A's age is in a particular range, but once you solve it, you see that his age makes sense. In this generalized puzzle I am giving you, you might be surprised how long wizards can live and keep their fertility.

Another difference with the original puzzle is that it is difficult to solve this one without a computer.

Once again the solution is unique. The number of the bus is 26. All buses with numbers below 25 do not give wizard A the opportunity to say ``no." If the bus is 26, the wizard is 3456 years old and has 7 children and 124 dolls. The ages of the children are one of the following:

\begin{itemize}
\item 1, 3, 3, 3, 4, 4, 8
\item 2, 2, 2, 2, 6, 6, 6
\end{itemize}

If the bus number is 27, we can re-use the solution for bus 26 and add one more child, aged 1. Thus wizard A could be 3456 years of age, could have 8 children, and could have 125 dolls. This bus number, however, has another solution for the younger age of 2560, 6 children and 165 dolls:

\begin{itemize}
\item 1, 4, 4, 4, 4, 10
\item 2, 2, 2, 5, 8, 8
\end{itemize}

You can try to continue to the next step of generalization and create another puzzle by adding the next symmetric polynomial on the ages of the children, for example, the sum of cubes. In this case, I do not know if the puzzle works: that is, if there is a unique bus number. I did not have the patience or the computer wisdom to do this. I conjecture that uniqueness will stop. Luck will run out, and we will not have a clean puzzle any more.

\end{document}